\documentclass[12pt]{article}
\usepackage{amssymb, amsfonts, amsmath, amsthm, url, graphicx}

\def\3{\subset }
\def\4{\subseteq }

\def\0{\leqno}

\def\barr{\begin{array}}
\def\earr{\end{array}}

\def\Z{{\rlap{$\kern2pt{\rm Z}$}{\rm Z}\,}}

\title{\bf A note on a metric associated to certain finite groups}
\author{Marius T\u arn\u auceanu}
\date{June 27, 2015}

\begin{document}

\maketitle

\begin{abstract}
In this short note we introduce a new metric on certain finite groups. It leads to a class of groups for which the element orders satisfy an interesting inequality. This extends the class ${\rm CP}_2$ studied in our previous paper \cite{16}.
\end{abstract}
\bigskip

\noindent{\bf MSC (2010):} Primary 20D10, 20D20; Secondary 20D15,
20D25, 20E34.

\noindent{\bf Key words:} finite groups, element orders, CP-groups, metrics.
\bigskip

\section{Introduction}

In group theory, and especially in geometric group theory, several metrics on a finite group $G$ have been studied (see e.g. \cite{2,4,6}). These are important because they give a way to measure the distance between any two elements of $G$. A new metric on certain finite groups $G$ will be presented in the following.

Let $d:G\times G\longrightarrow\mathbb{N}$ be the function defined by $$d(x,y)=o(xy^{-1})-1, \forall\, x,y\in G,$$where $o(a)$ is the order of $a\in G$. Then $d$ is a metric on $G$ if and only if $$o(ab)<o(a)+o(b), \forall\, a,b\in G.\0(*)$$Denote by ${\rm CP}_3$ the class of finite groups $G$ which satisfy $(*)$.

Remark that $d$ becomes an ultrametric on $G$ if and only if $$o(ab)\leq{\rm max}\{o(a),o(b)\}, \forall\, a,b\in G,$$that is $G$ belongs to the class ${\rm CP}_2$ studied in \cite{16}. Consequently, ${\rm CP}_2$ is properly contained in ${\rm CP}_3$ (an example of a finite group in ${\rm CP}_3$ but not in ${\rm CP}_2$ is the symmetric group $S_3$). This implies that ${\rm CP}_3$ contains all abelian $p$-groups, as well as the quaternion group $Q_8$ or the alternating group $A_4$, but at first sight it is difficult to describe all finite groups in this class. Their study is the main goal of our note.

Most of our notation is standard and will not be repeated here. Basic notions and results on group theory can be found in \cite{7,11,12,15}.

\section{Main results}

First of all, we observe that ${\rm CP}_3$ is closed under subgroups. On the other hand, since the cyclic group $\mathbb{Z}_6$ does not belong to ${\rm CP}_3$ (it has two elements of orders 2 and 3 whose sum is of order 6), we infer that ${\rm CP}_3$ is not closed under direct products or extensions.

\bigskip\noindent{\bf Remarks.}
\begin{itemize}
\item[\rm 1.] We know that ${\rm CP}_2\subset{\rm CP}_3$. Then, by Remark 1 of \cite{16}, we are able to indicate other three classes of finite $p$-groups, more large as the class of abelian $p$-groups, that are contained in ${\rm CP}_3$: regular $p$-groups (see Theorem 3.14 of \cite{15}, II, page 47), $p$-groups whose subgroup lattices are modular (see Lemma 2.3.5 of \cite{13}), and powerful $p$-groups for $p$ odd (see the main theorem of \cite{17}).
\item[\rm 2.] $Q_8$ is the smallest nonabelian $p$-group contained in ${\rm CP}_3$, while the dihedral group $D_8$ is the smallest $p$-group not contained in ${\rm CP}_3$. Note that all quaternion groups $Q_{2^n}$, $n\geq 4$, as well as all dihedral groups $D_{2n}$, $n\geq 4$, does not belong to ${\rm CP}_3$.
\item[\rm 3.] The groups $S_n$ with $n\geq 4$ are not contained in ${\rm CP}_3$ (for this it is enough to observe that $S_4\not\in {\rm CP}_3$: there are $a=(12)(34),b=(13)\in S_4$ such that $o(ab)=4\nless o(a)+o(b)=4$). Similarly, the groups $A_n$ with $n\geq 5$ are also not contained in ${\rm CP}_3$.
\end{itemize}

The following theorem gives a connection between ${\rm CP}_3$ and the well-known class of CP-groups (see e.g. \cite{1,3,5,8,9,10,14,18}).

\bigskip\noindent{\bf Theorem 1.} {\it ${\rm CP}_3$ is properly contained in ${\rm CP}$.}

\bigskip\noindent{\bf Proof.}  Assume that a group $G$ in ${\rm CP}_3$ contains an element $x$ whose order is not a prime power. Then there are two powers $a$ and $b$ of $x$ such that $o(a)=p$ and $o(b)=q$, where $p$ and $q$ are distinct primes. Since $a$ and $b$ commute, we have $o(ab)=pq$. By $(*)$ it follows that
$$pq<p+q,$$a contradiction.

Clearly, $S_4$ is contained in ${\rm CP}$, but not in ${\rm CP}_3$. This shows that the inclusion ${\rm CP}_3\subset{\rm CP}$ is proper.
\hfill\rule{1,5mm}{1,5mm}
\bigskip

A similar argument as in the above proof leads to the following property of groups in ${\rm CP}_3$.

\bigskip\noindent{\bf Theorem 2.} {\it Any abelian subgroup of a group in ${\rm CP}_3$ is a $p$-group. In particular, an abelian group is contained in ${\rm CP}_3$ if and only if it is a $p$-group.}
\smallskip

Our next result proves that the intersections of ${\rm CP}_2$ and ${\rm CP}_3$ with the class of $p$-groups are the same.

\bigskip\noindent{\bf Theorem 3.} {\it A $p$-group $G$ is contained in ${\rm CP}_3$ if and only if it is contained in ${\rm CP}_2$.}

\bigskip\noindent{\bf Proof.}  Assume that $G$ belongs to ${\rm CP}_3$ and let $p^n$ be its order. We will prove that for every $i=0,1,...,n$ the set $G_i=\{x \in G \mid o(x)\leq p^i\}$ is a normal subgroup of $G$. Let $x,y\in G_i$. Then $o(x),o(y)\leq p^i$ and therefore $o(xy)< 2p^i$ by $(*)$. On the other hand, we know that $G$ belongs to ${\rm CP}$ and so $o(xy)=p^j$ for some non-negative integer $j$. Thus $p^j<2p^i$, which leads to $j\leq i$, i.e. $xy\in G_i$. This proves that $G_i$ is a subgroup of $G$. Moreover, $G_i$ is normal in $G$ because the order map is constant on each conjugacy class. Then Theorem A of \cite{16} implies that $G$ belongs to ${\rm CP}_2$, as desired.
\hfill\rule{1,5mm}{1,5mm}
\bigskip

By \cite{14} we know that only eight nonabelian finite simple CP-groups exist: ${\rm PSL}(2,q)$ for $q=4,7,8,9,17$, ${\rm PSL}(3,4)$, ${\rm Sz}(8)$, and ${\rm Sz}(32)$. All these groups are not contained in ${\rm CP}_3$, as shows the following theorem.

\bigskip\noindent{\bf Theorem 4.} {\it ${\rm CP}_3$ contains no nonabelian finite simple group.}

\bigskip\noindent{\bf Proof.} Since the product of any two elements of order 2 of a group in ${\rm CP}_3$ can have order at most 3, we infer that ${\rm PSL}(2,q)$ does not belong to ${\rm CP}_3$ whenever $q\geq 4$ (note that ${\rm PSL}(2,2)\cong S_3$ and ${\rm PSL}(2,3)\cong A_4$ belong to ${\rm CP}_3$). ${\rm PSL}(3,4)$ has a subgroup isomorphic to ${\rm PSL}(2,4)\cong A_5$, and consequently it also does not belong to ${\rm CP}_3$. The same conclusion is obtained for the Suzuki groups ${\rm Sz}(8)$ and ${\rm Sz}(32)$ because they contain a subgroup isomorphic to $D_{10}$.
\hfill\rule{1,5mm}{1,5mm}
\smallskip

Inspired by Theorem 4 and by Corollary E of \cite{16} we came up with the following conjecture.

\bigskip\noindent{\bf Conjecture 5.} {\it ${\rm CP}_3$ is properly contained in the class of finite solvable groups.}
\bigskip

Finally, we indicate three natural problems concerning the class of finite groups introduced in our paper.

\bigskip\noindent{\bf Problem 1.} Study whether ${\rm CP}_3$ is closed under homomorphic images.

\bigskip\noindent{\bf Problem 2.} Determine the intersection between ${\rm CP}_3$ and ${\rm CP}_1$.

\bigskip\noindent{\bf Problem 3.} Give a precise description of the structure of finite groups contained in ${\rm CP}_3$.

\vspace*{5ex}\small

\hfill
\begin{minipage}[t]{5cm}
Marius T\u arn\u auceanu \\
Faculty of  Mathematics \\
``Al.I. Cuza'' University \\
Ia\c si, Romania \\
e-mail: {\tt tarnauc@uaic.ro}
\end{minipage}

\end{document}